\documentclass[fleqn]{mat01}
\usepackage{times,mathtimy,amssymb,latexsym,amsbsy}
\begin{document}

\setcounter{page}{339}
\firstpage{339}

\newtheorem{theore}{Theorem}
\renewcommand\thetheore{\arabic{section}.\arabic{theore}}
\newtheorem{theor}[theore]{\bf Theorem}
\newtheorem{rem}[theore]{Remark}
\newtheorem{propo}[theore]{\rm PROPOSITION}
\newtheorem{lem}[theore]{Lemma}
\newtheorem{definit}[theore]{\rm DEFINITION}
\newtheorem{coro}[theore]{\rm COROLLARY}
\newtheorem{exampl}[theore]{Example}
\newtheorem{case}{Case}

\def\corol{\trivlist \item[\hskip \labelsep{COROLLARY.}]}
\def\noteproof{\trivlist \item[\hskip \labelsep{\it Note added in Proof.}]}

\renewcommand{\theequation}{\thesection\arabic{equation}}

\title{Fixed point theory for composite maps on almost dominating extension
spaces}

\markboth{Ravi P Agarwal, Jong Kyu Kim and Donal O'Regan}{Fixed point theory for composite maps}

\author{RAVI P AGARWAL$^{1}$, JONG KYU KIM$^{2}$ and DONAL O'REGAN$^{3}$}

\address{$^{1}$Department of Mathematical Sciences, Florida Institute of
Technology, Melbourne, Florida 32901-6975, USA\\
\noindent $^{2}$Department of Mathematics, Kyungnam University, Masan,
Kyungnam 631-701, Korea\\
\noindent $^{3}$Department of Mathematics, National University of
Ireland, Galway, Ireland\\
\noindent E-mail: agarwal@fit.edu}

\volume{115}

\mon{August}

\parts{3}

\pubyear{2005}

\Date{MS received 4 February 2004; revised 1 February 2005}

\begin{abstract}
New fixed point results are presented for
${\cal U}_c^{\kappa}(X,X)$ maps in extension type spaces.
\end{abstract}


\maketitle

\section{Introduction}

In this paper we present new fixed point results in extension type
spaces. In particular, we present results for compact upper
semicontinuous ${\cal U}_c^{\kappa}(X,X)$ maps where $X$ is almost
ES dominating or more generally almost Schauder admissible
dominating (these concepts will be defined in \S 2). The results
in this paper improve those in the literature (see
\cite{3,4,5,9,11} and the references therein). A continuation
theorem is also discussed when the maps are between topological
vector spaces.

For the remainder of this section we present some definitions and
known results which will be needed throughout this paper. Suppose
$X$ and $Y$ are topological spaces. Given a class ${\cal X}$ of
maps, ${\cal X}(X,Y)$ denotes the set of maps $F\hbox{\rm :}\ X \rightarrow
2^Y$ (nonempty subsets of $Y$) belonging to ${\cal X}$, and ${\cal
X}_c$ the set of finite compositions of maps in ${\cal X}$. We let
\begin{equation*}
{\cal F}({\cal X})=\left\{ Z\hbox{\rm :}\  \hbox{Fix}\ F \neq
\emptyset\quad\hbox{for all}\ F \in {\cal X}(Z,Z) \right\},
\end{equation*}
where Fix $F$ denotes the set of fixed points of $F$.

${\cal U}$ will be the class of maps \cite{11} with the following
properties:

\begin{enumerate}
\renewcommand\labelenumi{(\roman{enumi})}
\leftskip .35pc
\item ${\cal U}$ contains the class ${\cal C}$ of single
valued continuous functions;

\item each $F \in {\cal U}_c$ is upper semicontinuous and
compact valued; and

\item $B^n \in {\cal F}({\cal U}_c)$ for all $n\in
\{1,2,\ldots\}$; here $B^n=\{x \in {\bf R}^n\hbox{\rm :}\ \|x\| \leq 1\}$.
\end{enumerate}
${\cal U}_c^{\kappa}(X,Y)$ will consist of all maps $F\hbox{\rm :}\ X
\rightarrow 2^Y$ such that for each $F$ and each nonempty compact
subset $K$ of $X$ there exists a map $G \in {\cal U}_c(K,Y)$ such
that $G(x) \subseteq F(x)$ for all $x\in K$.

Recall \cite{9} that ${\cal U}_c^{\kappa}$ is closed under
compositions. We also discuss special examples of ${\cal
U}_c^{\kappa}$ maps. Let $X$ and $Y$ be subsets of Hausdorff
topological vector spaces $E_1$ and $E_2$ respectively. We will
consider maps $F\hbox{\rm :}\ X \rightarrow K(Y)$; here $K(Y)$ denotes the
family of nonempty compact subsets of $Y$. We say $F\hbox{\rm :}\ X \rightarrow
K(Y)$ is {\it Kakutani} if $F$ is upper semicontinuous with convex
values. A nonempty topological space is said to be acyclic if all
its reduced \u Cech homology groups over the rationals are
trivial. Now $F\hbox{\rm :}\ X \rightarrow K(Y)$ is {\it acyclic} if $F$ is
upper semicontinuous with acyclic values. $F\hbox{\rm :}\ X \rightarrow K(Y)$
is said to be an {\it O'Neill} map if $F$ is continuous and if the
values of $F$ consist of one or $m$ acyclic components (here $m$
is fixed).

Given two open neighborhoods $U$ and $V$ of the origins in $E_1$
and $E_2$ repectively, a $(U,V)$-approximate continuous selection
of $F\hbox{\rm :}\ X \rightarrow K(Y)$ is a continuous function $s\hbox{\rm :}\ X
\rightarrow Y$ satisfying
\begin{equation*}
s(x) \in \left(F \left[(x+U) \cap X\right] +V \right) \cap Y
\quad\hbox{for every}\ x\in X.
\end{equation*}
We say $F\hbox{\rm :}\ X \rightarrow K(Y)$ is {\it approximable} if it is upper
semicontinuous and if its restriction $F|_K$ to any compact subset
$K$ of $X$ admits a $(U,V)$-approximate continuous selection for
every open neighborhood $U$ and $V$ of the origins in $E_1$ and
$E_2$ respectively.

For our next definition let $X$ and $Y$ be metric spaces. A
continuous single-valued map $p\hbox{\rm :}\ Y \rightarrow X$ is called a
Vietoris map if the following two conditions are satisfied:

\begin{enumerate}
\renewcommand\labelenumi{(\roman{enumi})}
\leftskip .15pc
\item for each $x\in X$, the set $p^{-1}(x)$ is acyclic;

\item $p$ is a proper map i.e. for every compact $A \subseteq X$ we
have that $p^{-1}(A)$ is compact.
\end{enumerate}

\begin{definit}$\left.\right.$\vspace{.5pc}

{\rm  \noindent A multifunction $\phi\hbox{\rm :}\ X \rightarrow K(Y)$ is {\it
admissible} (strongly) in the sense of Gorniewicz, if $\phi\hbox{\rm :}\ X
\rightarrow K(Y)$ is upper semicontinuous, and if there exists a
metric space $Z$ and two continuous maps $p\hbox{\rm :}\ Z \rightarrow X$ and
$q\hbox{\rm :}\ Z \rightarrow Y$ such that

\begin{enumerate}
\renewcommand\labelenumi{(\roman{enumi})}
\leftskip .15pc
\item $p$ is a Vietoris map, and

\item $\phi(x)=q(p^{-1}(x))$ for any $x\in X$.
\end{enumerate}

It should be noted that $\phi$ upper semicontinuous is redundant
in Definition~1.1. Notice the Kakutani maps, the acyclic maps, the
O'Neill maps, the approximable maps and the maps admissible in the
sense of Gorniewicz are examples of ${\cal U}_c^{\kappa}$ maps.}
\end{definit}

\section{Fixed point theory}

We begin with a result which extends results in the literature
\cite{3,11} for ${\cal U}_c^{\kappa}$ maps. In this paper by a
space we mean a Hausdorff topological space. Let $Q$ be a class of
topological spaces. A space $Y$ is an {\it extension space} for
$Q$ (written $Y \in ES(Q)$) if $\forall X \in Q$, $\forall K
\subseteq X$ closed in $X$, any continuous function $f_{0}\hbox{\rm :}\ K
\rightarrow Y$ extends to a continuous function $f\hbox{\rm :}\ X \rightarrow
Y$.  The following result was established in \cite{4}.

\setcounter{theore}{0}
\begin{theor}[\!]
Let $X\in ES(\hbox{compact})$ and $F \in {\cal U}_c^{\kappa}(X,X)$ a
compact map. Then $F$ has a fixed point.
\end{theor}

\looseness -1 We begin with topological vector spaces, so let $E$ be a Hausdorff
topological vector space and $X \subseteq E$. Let $V$ be a
neighborhood of the origin $0$ in $E$. $X$ is said to be {\it ES
$V$\!-dominated} if there exists a space $X_V \in
ES(\hbox{compact})$ and two continuous functions $r_{V}\hbox{\rm :}\ X_V
\rightarrow X$, $s_{V}\hbox{\rm :}\ X \rightarrow X_V$ such that $x-r_V\,s_V(x)
\in V$ for all $x\in X$. $X$ is said to be {\it almost ES
dominated} if $X$ is ES $V$-dominated for every neighborhood $V$
of the origin $0$ in $E$.

Any space that is dominated by a normed space (or more generally a
complete metric topological space admissible in the sense of Klee)
is almost ES dominated.

\begin{theor}[\!]
Let $X$ be a subset of a Hausdorff topological vector space $E$.
Also assume $X$ is almost ES dominated and that $F \in {\cal
U}_c^{\kappa}(X,X)$ is a compact closed map. Then $F$ has a fixed
point.
\end{theor}

\begin{proof}
Let ${\cal N}$ be a fundamental system of neighborhoods of the
origin $0$ in $E$ and $V \in {\cal N}$. Let $K=\overline{F(X)}$.
Now there exists $x_V \in ES(\hbox{compact})$ and two continuous
functions $r_{V}\hbox{\rm :}\ X_V \rightarrow X$, $s_{V}\hbox{\rm :}\ X \rightarrow X_V$ such
that $x-r_V s_V(x) \in V$ for all $x\in X$.

Let us look at the map $G_V=s_V Fr_V$. Since ${\cal U}_c^{\kappa}$
is closed under compositions we have that $G_V \in {\cal
U}_c^{\kappa}(X_V,X_V)$ is a compact map. Now Theorem~2.1
guarantees that there exists $z_V \in X_V$ with $z_V \in s_V Fr_V
(z_V)$ i.e. $z_V=s_V(x_V)$ for some $x_V \in Fr_V (z_V)$. In
particular $r_V(z_V)=r_V s_V(x_V)$ so $x_V \in F r_V (z_V)= Fr_V
s_V(x_V)$.

Let $y_V=r_V s_V(x_V)$. Then $x_V \in F y_V$ and since $x-r_V
s_V(x) \in V$ for all $x\in X$ we have $x_V-y_V \in V$. Now since
$K=\overline{F(X)}$ is compact we may assume without loss of
generality that there exists a $x$ with $x_V \rightarrow x$. Also
since $x_V-y_V \in V$ we have $y_V \rightarrow x$. These together
with the facts that $F$ is closed and $x_V \in F(y_V)$ implies
$x\in F(x)$.\hfill $\Box$ \end{proof}

Let $Q$ be a class of topological spaces and let $Y$ be a subset
of a Hausdorff topological vector space. $Y$ is a {\it Klee
approximate extension space} for $Q$ (written $Y \in {\it KAES}(Q)$) if
for all neighborhoods $V$ of $0$, $\forall X \in Q$, $\forall K
\subseteq X$ closed in $X$, and any continuous function $f_{0}\hbox{\rm :}\ K
\rightarrow Y$, there exists a continuous function $f_{V}\hbox{\rm :}\ X
\rightarrow Y$ with $f_V(x)-f_0(x) \in V$ for all $x\in K$.

Let $E$ be a Hausdorff topological vector space and $X \subseteq
E$. We say $X$ is {\it KAES admissible} if for every compact
subset $K$ of $X$ and every neighborhood $V$ of $0$ there exists a
continuous function $h_{V}\hbox{\rm :}\ K \rightarrow X$ such that

\begin{enumerate}
\renewcommand\labelenumi{(\roman{enumi})}
\leftskip .15pc
\item $x-h_V(x) \in V$ for all $x\in K$ and

\item $h_V(K)$ is contained in a subset $C$ of $X$ with $C \in
{\it KAES}(\hbox{compact})$.
\end{enumerate}

The following result was established in \cite{3}.

\begin{theor}[\!]
Let $X$ be a subset of a Hausdorff topological vector space $E$.
Also assume $X$ is KAES admissible and that $F \in {\cal
U}_c^{\kappa}(X,X)$ is a upper semicontinuous compact map with
closed values. Then $F$ has a fixed point.
\end{theor}

Let $E$ be a Hausdorff topological vector space and $X \subseteq
E$. Let $V$ be a neighborhood of the origin $0$ in $E$. $X$ is
said to be {\it KAES $V$-dominated} if there exists a KAES
admissible space $X_V$ and two continuous functions $r_{V}\hbox{\rm :}\ X_V
\rightarrow X$, $s_{V}\hbox{\rm :}\ X \rightarrow X_V$ such that $x-r_V\,s_V(x)
\in V$ for all $x\in X$. $X$ is said to be {\it almost KAES
dominated} if $X$ is KAES $V$-dominated for every neighborhood $V$
of the origin $0$ in $E$.

Essentially the same reasoning as in Theorem~2.2 (except here we
use Theorem~2.3) yields the following result.

\begin{theor}[\!]
Let $X$ be a subset of a Hausdorff topological vector space $E$.
Also assume $X$ is almost KAES dominated and that $F \in {\cal
U}_c^{\kappa}(X,X)$ is a compact upper semicontinuous map with
closed values. Then $F$ has a fixed point.
\end{theor}

Let $X$ be a subset of a Hausdorff topological vector space. Then
$X$ is said to be {\it $q$-almost KAES dominated} if any nonempty
compact convex subset $\Omega$ of $X$ is almost KAES dominated.

\begin{theor}[\!]
Let $X$ be a closed convex $q$-almost KAES dominated subset of a
Hausdorff topological vector space with $x_0 \in X$. Suppose $F
\in {\cal U}_c^{\kappa}(X,X)$ is a upper semicontinuous map with
closed values and assume the following condition holds{\rm :}
\begin{equation}
\hbox{If}\ A \subseteq X\ \hbox{with}\ A =\overline{co}(\{x_0\}
\cup F(A)), \ \hbox{then}\ A\ \hbox{is compact.}
\end{equation}
Then $F$ has a fixed point.
\end{theor}

\begin{proof}
Consider ${\cal F}$ the family of all closed, convex subsets $C$
of $X$ with $x_0 \in C$ and $F(x) \subseteq C$ for all $x\in C$
and let $C_0= \cap_{C \in {\cal F}}\,C$. Now \cite{2} guarantees
that
\begin{equation}
C_0=\overline{co}\,(\{x_0\} \cup F(C_0)).
\end{equation}
Now $(2.1)$ guarantees that $C_0$ is compact and $(2.2)$ implies
$F(C_0)\subseteq C_0$. Also $C_0$ is almost KAES dominated. In
addition $F|_{C_0} \in {\cal U}_c^{\kappa}(C_0,C_0)$ is a compact
upper semicontinuous map with closed values so Theorem~2.4
guarantees that there exists $x_0 \in C$ with $x_0 \in Fx_0$.
\hfill $\Box$
\end{proof}

Next we present a continuation theorem. Let $E$ be a Hausdorff
topological vector space, $C$ a closed convex subset of $E$, $U$
an open subset of $C$ and $0\in U$. We would like to consider maps
$F\hbox{\rm :}\ \overline{U} \rightarrow K(C)$ which are upper semicontinuous
and either (a)~Kakutani; (b)~acyclic; (c)~O'Neill;
(d)~approximable; or (e)~admissible (strongly) in the sense of
Gorniewicz. Here $\overline{U}$ denotes the closure of $U$ in $C$.

To take care of all the above maps (and even more general types)
we define as follows.

\setcounter{theore}{0}
\begin{definit}$\left.\right.$\vspace{.5pc}

{\rm \noindent $F\in GA(\overline{U}, C)$ if $F\hbox{\rm :}\ \overline{U}
\rightarrow K(C)$ is upper semicontinuous and satisfies condition
$(C)$.

We assume condition $(C)$ as:
\begin{equation}
\begin{cases}
\hbox{for any map}\ F \in GA(\overline{U},C)\ \hbox{and any
continuous single-valued}\\
\hbox{map}\ \mu\hbox{\rm :}\ \overline{U} \rightarrow [0,1]\ \hbox{we have
that}\ \mu F\ \hbox{satisfies condition}\ (C).
\end{cases}
\end{equation}
Certainly if condition $(C)$ means (a), (b), (c), (d) or (e) above,
then $(2.3)$ holds.

We assume the map $F\hbox{\rm :}\ \overline{U} \rightarrow K(C)$ satisfies one
of the following conditions:

\begin{enumerate}
\renewcommand\labelenumi{(H\arabic{enumi})}
\leftskip .7pc
\item $F\hbox{\rm :}\ \overline{U} \rightarrow K(C)$ is compact;

\item if $D \subseteq \overline{U}$ and $D \subseteq
\overline{co}\,(\{0\} \cup F(D))$ then $\overline{D}$ is compact.
\end{enumerate}
Fix $i\in \{1,2\}$.}
\end{definit}

\begin{definit}$\left.\right.$\vspace{.5pc}

{\rm \noindent We say $F\in GA^i(\overline{U}, C)$ if $F \in
GA(\overline{U}, C)$ satisfies (Hi).}
\end{definit}

\begin{definit}$\left.\right.$\vspace{.5pc}

{\rm \noindent We say $F\in GA^i_{\partial U}(\overline{U}, C)$ if
$F \in GA^i(\overline{U}, C)$ with $x \notin F(x)$ for $x \in
\partial U$; here $\partial U$ denotes the boundary of $U$
in $C$.}
\end{definit}

\begin{definit}$\left.\right.$\vspace{.5pc}

{\rm \noindent A map $F \in GA^i_{\partial U}(\overline{U}, C)$ is
essential in $GA^i_{\partial U}(\overline{U}, C)$ if for every $G
\in GA^i_{\partial U}(\overline{U}, C)$ with $G|_{\partial
U}=F|_{\partial U}$ there exists $x\in U$ with $x \in G(x)$.}
\end{definit}

The following result was established in \cite{1}.

\setcounter{theore}{5}
\begin{theor}[\!]
Fix $i\in \{1,2\}$. Let $E$ be a Hausdorff topological vector space{\rm
,} $C$ a closed convex subset of $E${\rm ,} $U$ an open subset of
$C${\rm ,} $0\in U$ and assume $(2.3)$ holds. Suppose $F\in
GA^i(\overline{U},C)$ and assume the following two conditions are
satisfied{\rm :}
\begin{equation}
\hbox{the zero map is essential in}\ GA^i_{\partial
U}(\overline{U}, C)
\end{equation}
and
\begin{equation}
x \notin \lambda Fx\quad \hbox{for}\ x\in \partial U\quad
\hbox{and}\quad \lambda \in (0,1].
\end{equation}
Then $F$ is essential in $GA^i_{\partial U}(\overline{U}, C)$.
\end{theor}

Next we discuss Assumption $(2.4)$. The results presented were
motivated partly by \cite{10}.

\setcounter{theore}{0}
\begin{exampl}
{\rm Let $i=1$. Suppose condition $(C)$ in Definition~2.1 means
$F\hbox{\rm :}\ \overline{U} \rightarrow K(C)$ is either (a)~Kakutani;
(b)~acyclic; (c)~O'Neill; or (d)~approximable, and assume $C$ is
almost KAES dominated. Then $(2.3)$ and $(2.4)$ hold.

To see $(2.4)$ let $\theta \in GA^1_{\partial U}(\overline{U}, C)$
with $\theta|_{\partial U}=\{0\}$. We must show that there exists
$x\in U$ with $x \in \theta(x)$. Let
\begin{equation*}
J(x)=
\begin{cases}
\theta(x), &x\in \overline{U} \\
\{0\}, &\hbox{otherwise.}
\end{cases}
\end{equation*}
It is easy to see that $J\hbox{\rm :}\ C \rightarrow K(C)$ is (a)~Kakutani if
$\theta$ is Kakutani, (b)~acyclic if $\theta$ is acyclic,
(c)~O'Neill if $\theta$ is O'Neill or (d)~approximable \cite{10}
if $\theta$ is approximable. Clearly $J\hbox{\rm :}\ C \rightarrow K(C)$ is
compact. Now Theorem~2.4 guarantees that there exists $x\in C$
with $x \in J(x)$. If $x \notin U$ we have $x\in J(x)=\{0\}$,
which is a contradiction since $0\in U$. Thus $x\in U$ so $x\in
J(x)=\theta(x)$.}
\end{exampl}

\begin{exampl}
{\rm Let $i=2$. Suppose condition $(C)$ in Definition~2.1 means
$F\hbox{\rm :}\ \overline{U} \rightarrow K(C)$ is either (a)~Kakutani;
(b)~acyclic; (c)`O'Neill; or (d)~approximable, and assume $C$ is
closed, convex and $q$-almost KAES dominated. Also suppose
\begin{equation}
\overline{co}\,(K)\ \hbox{is compact for any compact subset}\ K\
\hbox{of}\ C.
\end{equation}
Then $(2.3)$ and $(2.4)$ hold.

To see $(2.4)$ let $\theta \in GA^2_{\partial U}(\overline{U}, C)$
with $\theta|_{\partial U}=\{0\}$ and let $J$ be as in Example
2.1. Again $J\hbox{\rm :}\ C \rightarrow K(C)$ is (a)~Kakutani if $\theta$ is
Kakutani, (b)~acyclic if $\theta$ is acyclic, (c)~O'Neill if
$\theta$ is O'Neill, or (d)~approximable if $\theta$ is
approximable. Next we show that $J$ satisfies $(2.1)$ (with $F$
replaced by $J$ and $x_0$ by $0$). To see this let $D \subseteq C$
with $D=\overline{co}\,\left(\{0\} \cup J(D)\right)$. Then
\begin{equation}
D \subseteq \overline{co}\,\left(\{0\} \cup \theta(D \cap U)
\right),
\end{equation}
and so
\begin{equation*}
D \cap U \subseteq \overline{co}\,\left(\{0\} \cup \theta(D \cap U) \right).
\end{equation*}
Since $\theta\in GA^2(\overline{U},C)$ we have that $\overline{D
\cap U}$ is compact. In addition $\theta$ upper semicontinuous
guarantees that $\theta(\overline{D \cap U})$ is compact, and this
together with $(2.6)$ implies that $\overline{co}\,\left(\{0\}
\cup \theta(\overline{D \cap U}) \right)$ is compact. Now $(2.7)$
implies $\overline{D}\,(=D)$ is compact, so $(2.1)$ holds.
Theorem~2.5 guarantees that there exists $x\in C$ with $x \in
J(x)$. As in Example~2.1 we have $x\in U$ so $x\in
J(x)=\theta(x)$.}
\end{exampl}

Next we extend the results of this section to Hausdorff
topological spaces. First we gather together some preliminaries.
For a subset $K$ of a topological space $X$, we denote by
$\hbox{Cov}_{X}(K)$ the set of all coverings of $K$ by open sets
of $X$ (usually we write $\hbox{Cov}(K)=\hbox{Cov}_{X}(K)$). Given
a map $F\hbox{\rm :}\ X \rightarrow 2^X$ and $\alpha \in \hbox{Cov}(X)$, a
point $x\in X$ is said to be an $\alpha$-fixed point of $F$ if
there exists a member $U\in \alpha$ such that $x\in U$ and $F(x)
\cap U \neq \emptyset$. Given two maps $F,G\hbox{\rm :}\ X \rightarrow 2^Y$ and
$\alpha\in \hbox{Cov}(Y)$, $F$ and $G$ are said to be
$\alpha$-close if for any $x\in X$ there exists $U_x \in \alpha$,
$y\in F(x) \cap U_x$ and $w\in G(x) \cap U_x$.

The following result can be found in \cite{5}.

\setcounter{theore}{6}
\begin{theor}[\!]
Let $X$ be a regular topological space and $F\hbox{\rm :}\ X \rightarrow
2^X$ an upper semicontinuous map with closed values. Suppose there
exists a cofinal family of coverings $\theta \subseteq \hbox{\rm
Cov}_X(\overline{F(X)})$ such that $F$ has an $\alpha$-fixed point for
every $\alpha \in \theta$. Then $F$ has a fixed point.
\end{theor}

\setcounter{theore}{0}
\begin{rem}
{\rm From Theorem~2.7 in proving the existence of fixed points in
uniform spaces for upper semicontinuous compact maps with closed
values it suffices (p.~298 of \cite{6}) to prove the existence of
approximate fixed points (since open covers of a compact set $A$
admit refinements of the form $\{U[x]\hbox{\rm :}\  x\in A\}$ where $U$ is a
member of the uniformity (p.~199 of \cite{8}) and such refinements
form a cofinal family of open covers). Note also that uniform
spaces are regular (in fact completely regular) (p.~431 of \cite{7}).
Note in Theorem~2.7 if $F$ is compact valued then the
assumption that $X$ is regular can be removed. For convenience, in
this paper we will apply Theorem~2.7 only when the space is
uniform.}
\end{rem}

Let $X$ be a Hausdorff topological space and let $\alpha \in
\hbox{Cov}(X)$. $X$ is said to be {\it ES $\alpha$-dominated} if
there exists a space $X_{\alpha} \in ES(\hbox{compact})$ and two
continuous functions $r_{\alpha}\hbox{\rm :}\ X_{\alpha} \rightarrow X$,
$s_{\alpha}\hbox{\rm :}\ X \rightarrow X_{\alpha}$ such that
$r_{\alpha}s_{\alpha}\hbox{\rm :}\ X \rightarrow X$ and $i\hbox{\rm :}\ X \rightarrow X$
are $\alpha$-close. $X$ is said to be {\it almost ES dominated} if
$X$ is ES $\alpha$-dominated for each $\alpha \in
\hbox{Cov}(X)$.

\setcounter{theore}{7}
\begin{theor}[\!]
Let $X$ be a uniform space and let $X$ be almost ES dominated.
Also suppose $F \in {\cal U}_c^{\kappa}(X,X)$ is a compact upper
semicontinuous map with closed values. Then $F$ has a fixed point.
\end{theor}

\begin{proof}
Let $K=\overline{F(X)}$ and $\alpha\in \hbox{Cov}_X(K)$. Now
there exists $X_{\alpha} \in ES(\hbox{compact})$ and two
continuous functions $r_{\alpha}\hbox{\rm :}\ X_{\alpha} \rightarrow X$,
$s_{\alpha}\hbox{\rm :}\ X \rightarrow X_{\alpha}$ such that
$r_{\alpha}s_{\alpha}$ and $i$ are $\alpha$-close. Notice
$s_{\alpha}Fr_{\alpha} \in {\cal
U}_c^{\kappa}(X_{\alpha},X_{\alpha})$ so Theorem~2.1 guarantees
that there exists $z_{\alpha} \in X_{\alpha}$ with $z_{\alpha} \in
s_{\alpha}Fr_{\alpha}(z_{\alpha})$ i.e. $z_{\alpha}=
s_{\alpha}(x_{\alpha})$ for some $x_{\alpha} \in
Fr_{\alpha}(z_{\alpha})$. Notice $x_{\alpha} \in
Fr_{\alpha}s_{\alpha}(x_{\alpha})$ also. Let
$y_{\alpha}=r_{\alpha}s_{\alpha}(x_{\alpha})$ and note
$x_{\alpha} \in Fy_{\alpha}$. Since $r_{\alpha}s_{\alpha}$ and
$i$ are $\alpha$-close there exists $U_{\alpha} \in \alpha$ with
$r_{\alpha}s_{\alpha}(x_{\alpha}) \in U_{\alpha}$ and
$x_{\alpha} \in U_{\alpha}$ i.e. $y_{\alpha} \in U_{\alpha}$ and
$x_{\alpha} \in U_{\alpha}$. As a result $y_{\alpha} \in
U_{\alpha}$ and $F(y_{\alpha}) \cap U_{\alpha} \neq \emptyset$
since $x_{\alpha} \in F(y_{\alpha})$ and $x_{\alpha} \in
U_{\alpha}$. Thus $F$ has an $\alpha$-fixed point. The result now
follows from Theorem~2.7 (with Remark~2.1). \hfill $\Box$
\end{proof}

Let $Q$ be a class of topological spaces and $Y$ a subset of a
Hausdorff topological space. A space $Y$ is an {\it approximate
extension space} for $Q$ (written $Y \in AES(Q)$) if
$\forall\,\alpha \in \hbox{Cov}(Y)$, $\forall X \in Q$, $\forall
K \subseteq X$ closed in $X$, and any continuous function $f_0\hbox{\rm :}\ K
\rightarrow Y$, there exists a continuous function $f\hbox{\rm :}\ X
\rightarrow Y$ such that $f|_{K}$ is $\alpha$-close to $f_0$.

Let $X$ be a uniform space. Then $X$ is {\it Schauder admissible}
if for every compact subset $K$ of $X$ and every covering $\alpha
\in \hbox{Cov}_X(K)$, there exists a continuous function (called
the Schauder projection) $\pi_{\alpha}\hbox{\rm :}\ K \rightarrow X$ such that

\begin{enumerate}
\renewcommand\labelenumi{(\roman{enumi})}
\leftskip .15pc
\item $\pi_{\alpha}$ and $i\hbox{\rm :}\ K \hookrightarrow X$ are
$\alpha$-close;

\item $\pi_{\alpha}(K)$ is contained in a subset $C \subseteq X$
with $C \in AES(\hbox{compact})$. \end{enumerate}

The following result was established in \cite{9}.

\begin{theor}[\!]
Let $X$ be a uniform space and assume $X$ is Schauder admissible.
Suppose $F \in {\cal U}_c^{\kappa}(X,X)$ is a compact upper
semicontinuous map with closed values. Then $F$ has a fixed point.
\end{theor}

Let $X$ be a Hausdorff topological space and let $\alpha \in
\hbox{Cov}(X)$. $X$ is said to be {\it Schauder admissible
$\alpha$-dominated} if there exists a Schauder admissible space
$X_{\alpha}$ and two continuous functions $r_{\alpha}\hbox{\rm :}\ X_{\alpha}
\rightarrow X$, $s_{\alpha}\hbox{\rm :}\ X \rightarrow X_{\alpha}$ such that
$r_{\alpha}s_{\alpha}\hbox{\rm :}\ X \rightarrow X$ and $i\hbox{\rm :}\ X \rightarrow X$
are $\alpha$-close. $X$ is said to be {\it almost Schauder
admissible dominated} if $X$ is Schauder admissible
$\alpha$-dominated for each $\alpha \in \hbox{Cov}(X)$.

Essentially the same reasoning as in Theorem~2.8 (except here we
use Theorem~2.9) yields the following result.

\begin{theor}[\!]
Let $X$ be a uniform space and let $X$ be almost Schauder admissible
dominated. Also suppose $F \in {\cal U}_c^{\kappa}(X,X)$ is a compact
upper semicontinuous map with closed values. Then $F$ has a fixed point.
\end{theor}


\begin{thebibliography}{99}
\bibitem{1} Agarwal R P and O'Regan D, Homotopy and Leray--Schauder
principles for multimaps, {\it Nonlinear Analysis Forum} {\bf 7}
(2002) 103--111

\bibitem{2} Agarwal R P and O'Regan D, Collectively fixed point
theorems, {\it Nonlinear Analysis Forum} {\bf 7} (2002) 167--179

\bibitem{3} Agarwal R P and O'Regan D, Fixed point theory for
multimaps defined on admissible subsets of topological vector
spaces, {\it Comm. Math. (Prace Mat.)} {\bf 44} (2004) 1--10

\bibitem{4} Agarwal R P, O'Regan D and Park S, Fixed point theory for
multimaps in extension type spaces, {\it J. Korean Math. Soc.}
{\bf 39} (2002) 579--591

\bibitem{5} Ben-El-Mechaiekh H, The coincidence problem for compositions
of set valued maps, {\it Bull. Austral. Math. Soc.} {\bf 41}
(1990) 421--434

\bibitem{6} Ben-El-Mechaiekh H, Spaces and maps approximation and
fixed points, {\it J. Comput. Appl. Math.} {\bf 113}
(2000) 283--308

\bibitem{7} Engelking R, General Topology (Berlin: Heldermann Verlag)
(1989)

\bibitem{8} Kelley J L , General Topology (New York: D. Van Nostrand
Reinhold Co.) (1955)

\bibitem{9} O'Regan D, Fixed point theory on extension type spaces
and essential maps on topological spaces, {\it Fixed Point Theory
and Applications} {\bf 2004(1)} (2004) 13--20

\bibitem{10} O'Regan D, Admissible spaces and $\alpha$-fixed points,
{\it Nonlinear Funct. Anal. Appl.} (to appear)

\bibitem{11} Park S, A unified fixed point theory of multimaps on topological vector
spaces, {\it J. Korean Math. Soc.} {\bf 35} (1998) 803--829
\end{thebibliography}
\end{document}